\documentclass[11pt]{article}
\usepackage{amsfonts}
\usepackage{amsmath}
\usepackage{hyperref}

\input{tcilatex}
\begin{document}

\title{\textbf{On Necessary and Sufficient Conditions for Near-Optimal
Singular Stochastic Controls}}
\author{\textbf{Mokhtar Hafayed}\thanks{%
E-mail address: hafa.mokh@yahoo.com} \\
\emph{Laboratory of Applied Mathematics, Mohamed Khider University, }\\
\emph{\ \ Po Box 145, Biskra 07000, Algeria.} \and \textbf{Syed Abbas}%
\thanks{%
E-mail address: sabbas.iitk@gmail.com} \\
\emph{School of Basic Sciences, Indian Institute of Technology Mandi,}\\
\emph{\ Mandi-175001, India.} \and \textbf{Petr Veverka}\thanks{%
E-mail address: panveverka@seznam.cz} \\
\emph{Department of Mathematics, Faculty of Nuclear Sciences and}\\
\emph{\ \ Physical Engineering, Czech Technical University, Trojanova 13,}\\
\emph{\ \ Prague 120 00, EU-Czech Republic.}}
\maketitle

\begin{abstract}
This paper is concerned with\ necessary and sufficient conditions for
near-optimal singular stochastic controls for systems driven by a nonlinear
stochastic differential equations (SDEs in short). The proof of our result
is based on Ekeland's variational principle and some dilecate estimates of
the state and adjoint processes. This result is a generalization of Zhou's
stochastic maximum principle for near-optimaity to singular control problem.


\end{abstract}

\section{\textbf{Introduction}}

\noindent In this paper, we consider the singular stochastic control problem
for systems governed by nonlinear controlled diffusion of the type%
\begin{equation}
\left\{ 
\begin{array}{c}
dx_{t}=f\left( t,x_{t},u_{t}\right) dt+\sigma \left( t,x_{t},u_{t}\right)
dW_{t}+G_{t}d\eta _{t}, \\ 
\multicolumn{1}{l}{x_{s}=y,}%
\end{array}%
\right.   \tag{1.1}  \label{1.1}
\end{equation}%
where $(W_{t})_{t}$ is a standard $l-$dimentional Brownian motion defined on
the filtered probability space $(\Omega ,\tciFourier ,\left( \tciFourier
_{t}\right) _{t},\mathbb{P}).$ The minimized criteria associated with the
state equation (\ref{1.1}) is defined by 
\begin{equation}
J\left( s,y,u,\eta \right) =\mathbb{E}\left[ h\left( x_{T}\right)
+\int_{s}^{T}\ell \left( t,x_{t},u_{t}\right) dt+\int_{s}^{T}k_{t}d\eta _{t}%
\right] ,  \tag{1.2}  \label{1.2}
\end{equation}%
where $\mathbb{E}$ denotes the expectation with respect to $\mathbb{P}$, and
the value function is defined as%
\begin{equation}
V\left( s,y\right) =\inf_{\left( u,\eta \right) \in \mathbb{U}\left( \left[
s,T\right] \right) }\left\{ J\left( s,y,u,\eta \right) \right\} .  \tag{1.3}
\label{1.3}
\end{equation}%
\noindent The kind of stochastic control problem has been investigated
extensively, both by the Bellman's dynamic programming \ method \cite{belman}
and by Pontryagin's maximum principle \cite{pontry}. In this paper, we are
concerned by th second method. Peng \cite{peng} introduced the second-order
adjoint equation and obtained the global maximum principle of optimality, in
which the control is present in the both drift and diffusion
coefficients.Studying near-optimal controls makes a good sense as studying
optimal controls from both theoretical as well as applications point of
view. Many more near-optimal controls are available than optimal ones,
indeed, optimal controls my not even exist in many situations, while
near-optimal controls always exist. The near-optimal deterministic controls
problem has been treated by many authors, including \cite{ekeland, zhoo1,
zhoo2, gabasov}. Recently,\ in an interesting paper, Zhou \cite{zhoo}
established the second-order necessary as well as sufficient conditions for
near-optimal stochastic controls for general controlled diffusion with two
adjoint processes. The near-optimal control problem for systems descripted
by volterra integral equations has been studied by \cite{pan}. However,
Chighoub et al., \cite{chighoub} extended Zhou's maximum principle of
near-optimality to SDEs with jumps. The similar problem for systems driven
by forward-backward stochastic differential equations has been solved in 
\cite{bahlali}. For justification of establishing a theory of near-optimal
controls see (\cite{zhoo, zhoo1}, Introduction).

Singular stochastic control problem is an important and challenging class of
problems in control theory, it appear in various fields like mathematical
finance, problem of optimal consumption etc. Stochastic maximum principle
for singular controls was considered by\ many authors, see for instance \cite%
{alvarez1, alvarez2, cadenillas, dufour, dufour1, sbahlali1, haussmann2,
karatzas1}. The first version of maximum principle for singular stochastic
control problems was obtained by Cadenillas et al., \cite{cadenillas}. The
first-order weak stochastic maximum principle has been studied in \cite%
{sbahlali1}. In \cite{dufour} the authors derived stochastic maximum
principle where the singular part has a linear form. Sufficient conditions
for existence of optimal singular control have been obtained in \cite%
{dufour1}.

The main objective of this paper is to establish necessary as well as
sufficient conditions for near-optimal singular control for SDEs. The
control domain is not necessarily convex. The proof of our result is based
on Ekeland's variational principle \cite{ekeland}, and some delicate
estimates of the state and adjoint processes. Finally, as an illustration an
example is solved explicitly. This result permit us to extend Zhou's maximum
principle of near-optimality to singular control problem.

\noindent The paper is organized as follows. The assumptions and statment
of\ the control problem is given in the second section.\ In the third and
forth section, we establish the main result of this paper.

\section{\textbf{Assumptions and statement of the problem}}

\noindent We consider stochastic optimal control of the following kind. Let $%
T$ be a fixed strictly positive real number and $(\Omega ,\tciFourier
,\left\{ \tciFourier _{t}\right\} _{t},\mathbb{P})$ be a fixed filtered
probability space satisfying the usual conditions in which a $l-$dimentional
Brownian motion $W=\left\{ W_{t}:s\leq t\leq T\right\} $ with $s\in \left[
0,T\right] $ and $W_{s}=0$\ is defined. Let $\mathbb{A}_{1}$ be a closed
convex subset of $\mathbb{R}^{m}$ and $\mathbb{A}_{2}:=\left( \left[
0,\infty \right) \right) ^{m}.$ Let $\mathbb{U}_{1}$ be the class of
measurable adapted processes $u:\left[ s,T\right] \times \Omega \rightarrow 
\mathbb{A}_{1}$ and $\mathbb{U}_{2}$ is the class of measurable adapted
processes $\eta :\left[ 0,T\right] \times \Omega \rightarrow \mathbb{A}_{2}$.

\noindent \textbf{Definition 1. }\textit{An admissible control is a pair }$%
\left( u,\eta \right) $ \textit{of measurable} $\mathbb{A}_{1}\times \mathbb{%
A}_{2}$-valued, $\tciFourier _{t}-$\textit{adapted processes, such that}

\noindent \textit{1) }$\eta $\textit{\ is of bounded variation,
nondecreasing\ continuous on the left with right limits and }$\eta _{s}=0.$

\noindent \textit{2) }$\mathbb{E}\left[ \sup_{t\in \left[ s,T\right]
}\left\vert u_{t}\right\vert ^{2}+\left\vert \eta _{T}\right\vert ^{2}\right]
<\infty .$

\noindent We denote $\mathbb{U}=\mathbb{U}_{1}\times \mathbb{U}_{2},$ the
set of all admissible controls. Since $d\eta _{t}$ may be singular with
respect to Lebesgue measure $dt,$ we call $\eta $ the singular part of the
control and the process $u$ its absolutely continuous part.

\noindent Throughout this paper, we also assume that

\begin{description}
\item[(H1)] $f:\left[ 0,T\right] \times \mathbb{R}^{n}\mathbb{\times \mathbb{%
A}}_{1}\mathbb{\rightarrow R}^{n},$ $\sigma :\left[ 0,T\right] \times 
\mathbb{R}^{n}\mathbb{\times A}_{1}\mathbb{\rightarrow }\mathcal{M}_{n\times
l}\left( \mathbb{R}\right) $ and $\ell :\left[ 0,T\right] \times \mathbb{R}%
^{n}\mathbb{\times \mathbb{A}}_{1}\rightarrow \mathbb{R},$ are measurable in 
$(t,x,u,)$ and twice continuously differentiable in $x$, and there exists a
constant $C>0$ such that, for $\varphi =f,\sigma ,\ell :$%
\begin{equation}
\left\vert \varphi (t,x,u,)-\varphi (t,x^{\prime },u,)\right\vert
+\left\vert \varphi _{x}(t,x,u,)-\varphi _{x}(t,x^{\prime },u,)\right\vert
\leq C\left\vert x-x^{\prime }\right\vert .  \tag{2.1}  \label{2.1}
\end{equation}%
\begin{equation}
\left\vert \varphi (t,x,u,)\right\vert \leq C\left( 1+\left\vert
x\right\vert \right) .  \tag{2.2}  \label{2.2}
\end{equation}

\item[(H2)] $h:\mathbb{R}^{n}\mathbb{\rightarrow R}$ is twice continuously
differentiable in $x$, and there exists a constant $C>0$ such that%
\begin{equation}
\left\vert h(x)-h(x^{\prime }))\right\vert +\left\vert
h_{x}(x)-h_{x}(x^{\prime }))\right\vert \leq C\left\vert x-x^{\prime
}\right\vert .  \tag{2.3}  \label{2.3}
\end{equation}%
\begin{equation}
\left\vert h(x)\right\vert \leq C\left( 1+\left\vert x\right\vert \right) . 
\tag{2.4}  \label{2.4}
\end{equation}

\item[(H3)] $G:\left[ 0,T\right] \rightarrow \mathcal{M}_{n\times m}\left( 
\mathbb{R}\right) ,$ $k:\left[ 0,T\right] \rightarrow \left( \left[ 0,\infty
\right) \right) ^{m},$ for each $t\in \lbrack 0,T]:G$ is continuous and
bounded, and $k$ is continuous.
\end{description}

\noindent Under the above assumptions, the SDE (\ref{1.1}) has a unique
strong solution $x_{t}$ which is given by%
\begin{equation*}
x_{t}=y+\int_{s}^{t}f\left( r,x_{r},u_{r}\right) dr+\int_{s}^{t}\sigma
\left( r,x_{r},u_{r}\right) dW_{r}+\int_{s}^{t}G_{r}d\eta _{r},
\end{equation*}%
and by standard arguments it is easy to show that for any $q>0,$ it hold that

\begin{equation*}
\mathbb{E(}\sup_{t\in \left[ s,T\right] }\left\vert x_{t}\right\vert
^{q})<C\left( q\right) ,
\end{equation*}%
where $C\left( q\right) $ is a constant depending only on $q$ and the
functional $J$ is well defined.

\noindent For any $\left( u,\eta \right) \in \mathbb{U}$ and the
corresponding state trajectory $x$, we define the first-order adjoint
process $\Psi _{t}$ and the second-order adjoint process $Q_{t}$ as the ones
satisfying the following two backward SDEs respectively

\begin{equation}
\left\{ 
\begin{array}{l}
d\Psi _{t}=-\left[ f_{x}^{\ast }\left( t,x_{t},u_{t}\right) \Psi _{t}+\sigma
_{x}^{\ast }\left( t,x_{t},u_{t}\right) K_{t}+\ell _{x}\left(
t,x_{t},u_{t}\right) \right] dt\smallskip \\ 
\text{ \ \ \ \ \ \ \ }+K_{t}dW_{t}, \\ 
\Psi _{T}=h_{x}\left( x_{T}\right) ,%
\end{array}%
\right.  \tag{2.5}  \label{2.5}
\end{equation}%
and%
\begin{equation}
\left\{ 
\begin{array}{l}
dQ_{t}=-\left[ f_{x}^{\ast }\left( t,x_{t},u_{t}\right)
Q_{t}+Q_{t}f_{x}^{\ast }\left( t,x_{t},u_{t}\right) +\sigma _{x}^{\ast
}\left( t,x_{t},u_{t}\right) Q_{t}\sigma _{x}^{\ast }\left(
t,x_{t},u_{t}\right) \right. \smallskip \\ 
\text{ \ \ \ \ \ \ \ \ \ }\left. +\sigma _{x}^{\ast }\left(
t,x_{t},u_{t}\right) R_{t}+R_{t}\sigma _{x}\left( t,x_{t},u_{t}\right)
+\Gamma _{t}\right] dt+R_{t}dW_{t}, \\ 
\Psi _{T}=h_{xx}\left( x_{T}\right) ,%
\end{array}%
\right.  \tag{2.6}  \label{2.6}
\end{equation}%
where%
\begin{equation*}
\Gamma _{t}=\ell _{xx}\left( t,x_{t},u_{t}\right)
+\tsum\limits_{i=1}^{n}\left( \Psi _{t}^{i}f_{xx}^{i}\left(
t,x_{t},u_{t}\right) +K_{t}^{i}\sigma _{xx}^{i}\left( t,x_{t},u_{t}\right)
\right) .
\end{equation*}%
As is well known, under conditions (H1), (H2) and (H3) the first-order
adjoint equation\ (\ref{2.5}) admits one and only one $\tciFourier -$adapted
solution pair $\left( \Psi ,K\right) \in \mathbb{L}_{\tciFourier }^{2}\left( %
\left[ 0,T\right] ,\mathbb{R}^{n}\right) \times \mathbb{L}_{\tciFourier
}^{2}\left( \left[ 0,T\right] ,\mathbb{R}^{n}\right) $ and the second-order
adjoint equation (\ref{2.6}) admits one and only one $\tciFourier -$adapted
solution pair $\left( Q,R\right) \in \mathbb{L}_{\tciFourier }^{2}\left( %
\left[ 0,T\right] ,\mathbb{R}^{n\times n}\right) \times \mathbb{L}%
_{\tciFourier }^{2}\left( \left[ 0,T\right] ,\mathbb{R}^{n\times n}\right) .$
Moreover, since $f_{x},\sigma _{x},\ell _{x}$ and $h_{x}$ are bounded then
we have the following estimate%
\begin{equation*}
\mathbb{E}\left[ \sup_{s\leq t\leq T}\left\vert \Psi _{t}\right\vert
^{2}+\int_{s}^{T}\left\vert K_{t}\right\vert ^{2}dt\right. +\left.
\sup_{s\leq t\leq T}\left\vert Q_{t}\right\vert ^{2}+\int_{s}^{T}\left\vert
R_{t}\right\vert ^{2}\right] \leq C.
\end{equation*}

\noindent Define the usual Hamiltonian%
\begin{equation}
H\left( t,x,u,p,q\right) :=-pf\left( t,x,u\right) -q\sigma \left(
t,x,u\right) -\ell \left( t,x,u\right) ,  \tag{2.7}  \label{2.7}
\end{equation}%
for $(t,x,u)\in \lbrack s,T]\times \mathbb{R}^{n}\times \mathbb{A}_{1}$.
Furthermore, we define the $\mathcal{H}$ functional corresponding to a given
admissible pair $\left( x,u\right) $ as\ follows%
\begin{eqnarray*}
\mathcal{H}^{\left( x,u\right) }(t,x,u) &=&H\left( t,x,u,\Psi
_{t},K_{t}-Q_{t}\sigma \left( t,x,u\right) \right) \\
&&-\frac{1}{2}\sigma ^{\ast }\left( t,x,u\right) Q_{t}\sigma \left(
t,x,u\right) ,
\end{eqnarray*}%
for $(t,x,u,p,q)\in \lbrack s,T]\times \mathbb{R}^{n}\times \mathbb{A}_{1}%
\mathbb{\times R}^{n}\times \mathbb{R}^{n}$, where $\Psi _{t},$ $K_{t}$ and $%
Q_{t}$ are determined by adjoint equations (\ref{2.5}) and (\ref{2.6})
corresponding to $(x,u).$

\noindent Befor concluding this section, let us recall Ekeland variational
principle and the Clarke's generalised gradient as follows

\noindent \textbf{Lemma 1. }(\textit{Ekeland's Lemma \cite{ekeland} })%
\textit{\ Let }$(F,$\textit{\ }$\rho )$\textit{\ be a complete metric space
and }$f:F\rightarrow \overline{\mathbb{R}}$\textit{\ be a lower
semi-continuous function which is bounded below. For a given }$\varepsilon
>0 $\textit{, suppose that }$u^{\varepsilon }\in F$\textit{\ satisfying }$%
f\left( u^{\varepsilon }\right) \leq \inf \left( f\right) +\varepsilon $%
\textit{, then for any }$\lambda >0$\textit{, there exists }$u^{\lambda }\in
F$\textit{\ such that}

\noindent \textit{1)} $f\left( u^{\lambda }\right) \leq f\left(
u^{\varepsilon }\right) .$

\noindent \textit{2) }$\rho \left( u^{\lambda },u^{\varepsilon }\right) \leq
\lambda .$

\noindent \textit{3)\ }$f\left( u^{\lambda }\right) \leq f\left( u\right) +%
\dfrac{\varepsilon }{\lambda }\rho \left( u,u^{\lambda }\right) $\textit{\
for all }$u\in F.$

\noindent For $u,v\in \mathbb{U}$ To apply Ekeland's variational principle%
\textit{\ }to our problem, we define a distance function $\rho $ on the
space of admissible controls such that $\left( \mathbb{U},d\right) $ becomes
a complete metric space. To achieve this goal, we define for any $\left(
u,\eta \right) \ $and $\left( v,\xi \right) \in \mathbb{U}:$%
\begin{equation}
d\left( \left( u,\eta \right) ,\left( v,\xi \right) \right) =d_{1}\left(
u,v\right) +d_{2}\left( \eta ,\xi \right) ,  \tag{2.8}  \label{2.8}
\end{equation}%
where%
\begin{equation}
d_{1}\left( u,v\right) =\mathbb{P\otimes }dt\left\{ \left( w,t\right) \in
\Omega \times \left[ 0,T\right] :u\left( w,t\right) \neq v\left( w,t\right)
\right\} ,  \tag{2.9}  \label{2.9}
\end{equation}%
and%
\begin{equation}
d_{2}\left( \eta ,\xi \right) =\left[ \mathbb{E}(\sup_{t\in \left[ s,T\right]
}\left\vert \eta _{t}-\xi _{t}\right\vert ^{2})\right] ^{\frac{1}{2}}, 
\tag{2.10}  \label{2.10}
\end{equation}%
here $\mathbb{P\otimes }dt$ is the product measure of $\mathbb{P}$ with the
Lebesgue measure $dt$ on $\left[ s,T\right] .$ It is easy to see that $%
\left( \mathbb{U}_{2},d_{2}\right) $ is a complete metric space. Moreover,
it has been shown in Yong \textit{et al.,} (\cite{yong} pp. 146-147) that $%
\left( \mathbb{U}_{1},d_{1}\right) $ is a complete metric space. Hence $%
\left( \mathbb{U},d\right) $ as a product of two complete metric spaces is a
complete metric space under $d$.

\noindent \textbf{Definition 1. (}Clarke's generalized gradient \cite{clarke}%
\textbf{)} Let $E$ be a convex set in $\mathbb{R}^{n}$ and let $%
f:E\rightarrow \mathbb{R}$ be a locally Lipschitz function. The Clarke's
generalized gradient of $f$ at $x\in E$, denoted by $\partial _{x}f$, is a
set defined by%
\begin{equation*}
\partial _{x}f\left( x\right) =\left\{ \xi \in \mathbb{R}^{n}:\left\langle
\xi ,v\right\rangle \leq \lim_{y\rightarrow x}\sup_{t\rightarrow 0}\frac{%
f\left( y+tv\right) -f\left( y\right) }{t}.\text{ }\forall v\in \mathbb{R}%
\text{ and }y,\text{ }\left( y+tv\right) \in E\right\} .
\end{equation*}

\section{Necessary conditions for near-optimal singular control}

\noindent Our purpose in this paper is to establish second-order necessary
and sufficient conditions for near-optimal singular control for systems
governed by nonlinear SDEs. It is worth montioning that optimal singular
controls may not even exist in many situations, while near-optimal singular
controls always exists. Ekeland's variational principle \cite{ekeland} is
applied to prove our maximum principle. The proof follows the general ideas
as in (\cite{zhoo1, zhoo2, zhoo}) where similar results are obtained for
other class of controls.

\noindent We give the definition of near-optimal control as given in Zhou (%
\cite{zhoo}, Definition 2.1 and Definition 2.2).

\noindent \textbf{Definition 2.} For a given $\varepsilon >0$\ the
admissible control $\left( u^{\varepsilon },\eta ^{\varepsilon }\right) $\
is near-optimal\ if%
\begin{equation}
\left\vert J\left( s,y,u^{\varepsilon },\eta ^{\varepsilon }\right) -V\left(
s,y\right) \right\vert \leq \mathcal{O}\left( \varepsilon \right) , 
\tag{3.1}  \label{3.1}
\end{equation}%
where $\mathcal{O}\left( .\right) $\ is a function of $\varepsilon $\
satisfying $\lim_{\varepsilon \rightarrow 0}\mathcal{O}\left( \varepsilon
\right) =0.$ The estimater $\mathcal{O}\left( \varepsilon \right) $ is
called an error bound. If $\mathcal{O}\left( \varepsilon \right)
=C\varepsilon ^{\delta }$ for some $\delta >0$ independent of the constant $%
C $ then $\left( u^{\varepsilon },\eta ^{\varepsilon }\right) $ is called
near-optimal control with order $\varepsilon ^{\delta }.$ If $\mathcal{O}%
\left( \varepsilon \right) =\varepsilon $ the admissible control $\left(
u^{\varepsilon },\eta ^{\varepsilon }\right) $ called $\varepsilon -$optimal.

\noindent Our first Lemma below, is deals with the continiuity of the state
processes under distance $\rho $

\noindent \textbf{Lemma 2. }\textit{If }$x_{t}^{u,\eta }$\textit{\ and }$%
x_{t}^{v,\eta }$\textit{\ be the solution of the state equation (\ref{1.1})
associated respectively with }$u$ and $v$.\textit{\ For any }$0<\alpha <1$
and $\beta >0$ satisfying $\alpha \beta <1$,\textit{\ there exists a
positive constants }$C=C\left( T,\alpha ,\beta \right) $\textit{\ such that}%
\begin{equation}
\mathbb{E(}\sup_{s\leq t\leq T}\left\vert x_{t}^{u,\eta }-x_{t}^{v,\xi
}\right\vert ^{2\beta })\leq Cd_{1}^{\alpha \beta }\left( u,v\right) . 
\tag{3.2}  \label{3.2}
\end{equation}

\noindent \textbf{Proof.} First, we assume that $\beta \geq 1$. Using
Burkholder-Davis-Gundy inequality for the martingale part, we can compute,
for any $r\geq s$%
\begin{eqnarray*}
\mathbb{E}\left[ \sup_{s\leq t\leq r}\left\vert x_{t}^{u,\eta }-x_{t}^{v,\xi
}\right\vert ^{2\beta }\right] &\leq &C\mathbb{E(}\int_{s}^{r}\left\{
\left\vert f\left( t,x_{t}^{u,\eta },u_{t}\right) -f\left( t,x_{t}^{v,\xi
},v_{t}\right) \right\vert ^{2\beta }\right. \\
&&\left. +\int_{s}^{r}\left\vert \sigma \left( t,x_{t}^{u,\eta
},u_{t}\right) -\sigma \left( t,x_{t}^{v,\xi },v_{t}\right) \right\vert
^{2\beta }\right\} dt \\
&&+C\mathbb{E}\left\vert \eta _{T}-\xi _{T}\right\vert ^{2\beta }, \\
&\leq &C\mathbb{E(}\int_{s}^{r}\left\{ \left\vert f\left( t,x_{t}^{u,\eta
},u_{t}\right) -f\left( t,x_{t}^{u,\eta },v_{t}\right) \right\vert ^{2\beta
}\right. \\
&&\left. +\int_{s}^{r}\left\vert \sigma \left( t,x_{t}^{u,\eta
},u_{t}\right) -\sigma \left( t,x_{t}^{u,\eta },v_{t}\right) \right\vert
^{2\beta }\right\} \mathbf{\chi }_{u_{t}\neq v_{t}}\left( t\right) dt \\
&&+C\mathbb{E(}\int_{s}^{r}\left\{ \left\vert f\left( t,x_{t}^{u,\eta
},v_{t}\right) -f\left( t,x_{t}^{u,\eta },v_{t}\right) \right\vert ^{2\beta
}\right. \\
&&\left. +\int_{0}^{t}\left\vert \sigma \left( t,x_{t}^{u,\eta
},v_{t}\right) -\sigma \left( t,x_{t}^{u,\eta },v_{t}\right) \right\vert
^{2\beta }\right\} \\
&&+C\mathbb{E}\left\vert \eta _{T}-\xi _{T}\right\vert ^{2\beta },
\end{eqnarray*}%
now arguing as in (\cite{zhoo}, Lemma 3.1) taking $b=\frac{1}{\alpha \beta }%
>1$ and $a>1$ such that $\frac{1}{a}+\frac{1}{b}=1,$ and applying
Cauchy-Schwarz inequality, we get

\begin{eqnarray*}
&&\mathbb{E}\int_{s}^{r}\left\vert f\left( t,x_{t}^{u,\eta },u_{t}\right)
-f\left( t,x_{t}^{u,\eta },v_{t}\right) \right\vert ^{2\beta }\mathbf{\chi }%
_{u_{t}\neq v_{t}}\left( t\right) dt \\
&\leq &\left\{ \mathbb{E}\int_{s}^{r}\left\vert f\left( t,x_{t}^{u,\eta
},u_{t}\right) -f\left( t,x_{t}^{u,\eta },v_{t}\right) \right\vert ^{2\beta
a}dt\right\} ^{\frac{1}{a}}\times \left\{ \mathbb{E}\int_{s}^{r}\mathbf{\chi 
}_{u_{t}\neq v_{t}}\left( t\right) dt\right\} ^{\frac{1}{b}},
\end{eqnarray*}%
using definition of $d_{1}$ and linear growth condition on $f$ we obtain

\begin{eqnarray*}
&&\mathbb{E}\int_{s}^{r}\left\vert f\left( t,x_{t}^{u,\eta },u_{t}\right)
-f\left( t,x_{t}^{u,\eta },v_{t}\right) \right\vert ^{2\beta }\mathbf{\chi }%
_{u_{t}\neq v_{t}}\left( t\right) dt \\
&\leq &C\left\{ \mathbb{E}\int_{s}^{r}\left( 1+\left\vert x_{t}^{u,\eta
}\right\vert ^{2\beta a}\right) dt\right\} ^{\frac{1}{a}}d_{1}\left(
u,v\right) ^{\alpha \beta } \\
&\leq &Cd_{1}\left( u,v\right) ^{\alpha \beta }.
\end{eqnarray*}%
Similarly, we can prove%
\begin{equation}
\mathbb{E}\int_{s}^{r}\left\vert \sigma \left( t,x_{t}^{u,\eta
},u_{t}\right) -\sigma \left( t,x_{t}^{u,\eta },v_{t}\right) \right\vert
^{2\beta }\mathbf{\chi }_{u_{t}\neq v_{t}}\left( t\right) dt\leq
Cd_{1}\left( u,v\right) ^{\alpha \beta }.  \tag{3.3}  \label{3.3}
\end{equation}%
Therefore, by using assumption (H1), we conclued that 
\begin{equation*}
\mathbb{E}(\underset{s\leq t\leq r}{\sup }\left\vert x_{t}^{u,\eta
}-x_{t}^{v,\xi }\right\vert ^{2\beta })\leq C\left\{ \mathbb{E}\int_{s}^{r}%
\underset{s\leq r\leq \theta }{\sup }\left\vert x_{t}^{u,\eta }-x_{t}^{v,\xi
}\right\vert ^{2\beta }d\theta +\mathbb{E}\left\vert \eta _{T}-\xi
_{T}\right\vert ^{2\beta }+d_{1}\left( u,v\right) ^{\alpha \beta }\right\} .
\end{equation*}%
Hence (\ref{3.1}) follows immediately from definition 1 and Gronwall's
inequality.

\noindent Now we assume $0\leq \beta <1$. Since $\frac{2}{\alpha }>1$ then
the Cauchy-Schwarz inequality yields%
\begin{eqnarray*}
\mathbb{E}(\underset{s\leq t\leq T}{\sup }\left\vert x_{t}^{u,\eta
}-x_{t}^{v,\xi }\right\vert ^{2\beta }) &\leq &\left[ \mathbb{E}(\underset{%
s\leq t\leq T}{\sup }\left\vert x_{t}^{u,\eta }-x_{t}^{v,\xi }\right\vert
^{2})\right] ^{\beta } \\
&\leq &\left[ Cd_{1}\left( u,v\right) ^{\alpha }\right] ^{\beta } \\
&\leq &Cd_{1}\left( u,v\right) ^{\alpha \beta }.
\end{eqnarray*}%
This completes the proof of Lemma 2.

\noindent \textbf{Lemma 3. }\textit{For any }$0<\alpha <1$\textit{\ and }$%
1<\beta <2$\textit{\ satisfying }$\left( 1+\alpha \right) \beta <2$\textit{,
there exist a positive constant }$C=C\left( \alpha ,\beta \right) $ such
that for any $\left( u,\eta \right) $\textit{, }$\left( v,\xi \right) \in 
\mathbb{U}\left( \left[ s,T\right] \right) $, along with the correspending
trajectories $x^{u,\eta }$, $x^{v,\xi }$ and the solutions\textit{\ }$\left(
\Psi ,K,Q,R\right) ,$\textit{\ }$\left( \Psi ^{\prime },K^{\prime
},Q^{\prime },R^{\prime }\right) $\textit{\ of the corresponding adjoint
equations, it holds that}

\begin{eqnarray}
\mathbb{E}\int_{0}^{T}\mathbb{(}\left\vert \Psi _{t}-\Psi _{t}^{\prime
}\right\vert ^{\beta }+\left\vert K_{t}-K_{t}^{\prime }\right\vert ^{\beta
})dt &\leq &Cd_{1}^{\frac{\alpha \beta }{2}}\left( u,v\right) .  \TCItag{3.4}
\label{3.4} \\
\mathbb{E}\int_{0}^{T}\mathbb{(}\left\vert Q_{t}-Q_{t}^{\prime }\right\vert
^{\beta }+\left\vert R_{t}-R_{t}^{\prime }\right\vert ^{\beta })dt &\leq
&Cd_{1}^{\frac{\alpha \beta }{2}}\left( u,v\right) .  \TCItag{3.5}
\label{3.5}
\end{eqnarray}

\noindent \textbf{Proof.} Since the adjoint processes are independant to
singular part, we use similar argument as in Zhou (\cite{zhoo} Lemma 3.2).

Now we are able to state and prove the necessary conditions for near-optimal
singular control for our problem, which is the main result in this paper.

\noindent Let $\left( \Psi ^{\varepsilon },K^{\varepsilon }\right) $\textit{%
\ and }$\left( Q^{\varepsilon },R^{\varepsilon }\right) $\textit{\ be the
solution of adjoint equations (\ref{2.5}) and (\ref{2.6}) respectively
corresponding to }$\left( x^{\varepsilon },\left( u^{\varepsilon },\eta
^{\varepsilon }\right) \right) .$

\noindent \textbf{Theorem 1.}\textit{\ (Maximum principle for any
near-optimal singular control).\ For any }$\delta \in (0,\frac{1}{3}],$%
\textit{\ and any near-optimal singular control }$\left( u^{\varepsilon
},\eta ^{\varepsilon }\right) $\textit{\ there exists a positive constant }$%
C=C\left( \delta \right) $\textit{\ such that for\ each }$\varepsilon >0$

\begin{equation}
\left\{ 
\begin{array}{c}
-C\varepsilon ^{\delta }\leq \mathbb{E}\int_{s}^{T}\left\{ \frac{1}{2}\left(
\sigma \left( t,x_{t}^{\varepsilon },u\right) -\sigma \left(
t,x_{t}^{\varepsilon },u_{t}^{\varepsilon }\right) \right) ^{\ast
}Q_{t}^{\varepsilon }\left( \sigma \left( t,x^{\varepsilon },u\right)
-\sigma \left( t,x_{t}^{\varepsilon },u_{t}^{\varepsilon }\right) \right)
\right. \smallskip \\ 
+\Psi _{t}^{\varepsilon }\left( f\left( t,x_{t}^{\varepsilon },u\right)
-f\left( t,x_{t}^{\varepsilon },u_{t}^{\varepsilon }\right) \right)
+K_{t}^{\varepsilon }\left( \sigma \left( t,x_{t}^{\varepsilon },u\right)
-\sigma \left( t,x_{t}^{\varepsilon },u_{t}^{\varepsilon }\right) \right)
\smallskip \\ 
\left. +\left( \ell \left( t,x_{t}^{\varepsilon },u\right) -\ell \left(
t,x_{t}^{\varepsilon },u_{t}^{\varepsilon }\right) \right) \right\} dt,%
\end{array}%
\right.  \tag{3.6}  \label{3.6}
\end{equation}%
\textit{and}%
\begin{equation}
-C\varepsilon ^{\delta }\leq \mathbb{E}\left[ \int_{0}^{T}(k_{t}+G_{t}^{\ast
}\Psi _{t}^{\varepsilon })d\left( \eta -\eta ^{\varepsilon }\right) _{t}%
\right] .  \tag{3.7}  \label{3.7}
\end{equation}

\noindent \textbf{Proof.} By using Ekeland's variational principle with $%
\lambda =\varepsilon ^{\frac{2}{3}},$ there is an admissible pair $\left( 
\overline{x}^{\varepsilon },\left( \overline{u}^{\varepsilon },\overline{%
\eta }^{\varepsilon }\right) \right) $ such that for any $\left( u,\eta
\right) \in \mathbb{U}:$%
\begin{equation}
\rho \left( \left( u^{\varepsilon },\eta ^{\varepsilon }\right) ,\left( 
\overline{u}^{\varepsilon },\overline{\eta }^{\varepsilon }\right) \right)
\leq \varepsilon ^{\frac{2}{3}}.  \tag{3.8}  \label{3.8}
\end{equation}%
and

\begin{equation*}
J\left( s,y,u^{\varepsilon },\eta ^{\varepsilon }\right) \leq J\left(
s,y,u^{\varepsilon },\eta ^{\varepsilon }\right) +\varepsilon ^{\frac{1}{2}%
}d\left( \left( u,\eta \right) ,\left( \overline{u}^{\varepsilon },\overline{%
\eta }^{\varepsilon }\right) \right) .
\end{equation*}%
Notice that $\left( u^{\varepsilon },\eta ^{\varepsilon }\right) $ which is
near-optimal for the initial cost $J$ defined in (\ref{1.2}) is optimal for
the new cost $J^{\varepsilon }$ given by%
\begin{equation}
J^{\varepsilon }\left( s,y,u,\eta \right) =J\left( s,y,u,\eta \right)
+\varepsilon ^{\frac{1}{3}}d\left( \left( u,\eta \right) ,\left( \overline{u}%
^{\varepsilon },\overline{\eta }^{\varepsilon }\right) \right) .  \tag{3.9}
\label{3.9}
\end{equation}%
then we have%
\begin{equation*}
J^{\varepsilon }\left( s,y,\overline{u}^{\varepsilon },\overline{\eta }%
^{\varepsilon }\right) \leq J^{\varepsilon }\left( s,y,u,\eta \right) \text{
for any }\left( u,\eta \right) \in \mathbb{U}\left( \left[ s,T\right]
\right) ,
\end{equation*}%
Next, we use the spike variation techniques for $\overline{u}^{\varepsilon }$%
to drive the first variational inequality and we use convex perturbation for 
$\overline{\eta }^{\varepsilon }$ as follows

\noindent \textbf{First variational inequality: }For any $\theta >0,$ we
define the following strong perturbation $(\overline{u}_{t}^{\varepsilon
,\theta },\overline{\eta }_{t}^{\varepsilon })\in \mathbb{U}:$

\begin{equation}
(\overline{u}^{\varepsilon ,\theta },\overline{\eta }_{t}^{\varepsilon
})=\left\{ 
\begin{array}{l}
\left( u,\overline{\eta }_{t}^{\varepsilon }\right) \text{, }t\in \left[
t_{0},t_{0}+\theta \right] ,\smallskip \\ 
\left( \overline{u}_{t}^{\varepsilon },\overline{\eta }_{t}^{\varepsilon
}\right) \text{, otherwise}.%
\end{array}%
\right.  \tag{3.10}  \label{3.10}
\end{equation}%
The fact that%
\begin{equation}
J^{\varepsilon }\left( s,y,\overline{u}^{\varepsilon },\overline{\eta }%
^{\varepsilon }\right) \leq J^{\varepsilon }(s,y,\overline{u}^{\varepsilon
,\theta },\overline{\eta }^{\varepsilon }),  \tag{3.11}  \label{3.11}
\end{equation}%
and%
\begin{equation*}
d((\overline{u}^{\varepsilon ,\theta },\overline{\eta }_{t}^{\varepsilon
}),\left( \overline{u}^{\varepsilon },\overline{\eta }_{t}^{\varepsilon
}\right) )=d_{1}(\overline{u}^{\varepsilon ,\theta },\overline{u}%
^{\varepsilon })\leq \theta ,
\end{equation*}%
imply that%
\begin{equation*}
J(s,y,\overline{u}^{\varepsilon ,\theta },\overline{\eta }^{\varepsilon
})-J\left( s,y,\overline{u}^{\varepsilon },\overline{\eta }^{\varepsilon
}\right) \geq -\theta \varepsilon ^{\frac{1}{3}}.
\end{equation*}%
Since the diference $J(s,y,\overline{u}^{\varepsilon ,\theta },\overline{%
\eta }^{\varepsilon }))-J\left( s,y,\overline{u}^{\varepsilon },\overline{%
\eta }^{\varepsilon }\right) $ is independant to the singular part, the
near-maximum condition (\ref{3.6}) follows by applying similar argument as
in Zhou (\cite{zhoo}), we get%
\begin{equation}
\begin{array}{l}
-C\varepsilon ^{\frac{1}{3}}\leq \mathbb{E}\int_{s}^{T}\left\{ \frac{1}{2}%
\left( \sigma \left( t,\overline{x}_{t}^{\varepsilon },u\right) -\sigma
\left( t,\overline{x}_{t}^{\varepsilon },\overline{u}_{t}^{\varepsilon
}\right) \right) ^{\ast }\overline{Q}_{t}^{\varepsilon }\left( \sigma \left(
t,\overline{x}^{\varepsilon },u\right) -\sigma \left( t,\overline{x}%
_{t}^{\varepsilon },\overline{u}_{t}^{\varepsilon }\right) \right) \right.
\smallskip \\ 
\text{ \ \ \ \ \ \ \ \ \ \ \ }+\overline{\Psi }_{t}^{\varepsilon }\left(
f\left( t,\overline{x}_{t}^{\varepsilon },u\right) -f\left( t,\overline{x}%
_{t}^{\varepsilon },\overline{u}_{t}^{\varepsilon }\right) \right) +%
\overline{K}_{t}^{\varepsilon }\left( \sigma \left( t,\overline{x}%
_{t}^{\varepsilon },u\right) -\sigma \left( t,\overline{x}_{t}^{\varepsilon
},\overline{u}_{t}^{\varepsilon }\right) \right) \smallskip \\ 
\text{ \ \ \ \ \ \ \ \ \ \ \ }\left. +\left( \ell \left( t,\overline{x}%
_{t}^{\varepsilon },u\right) -\ell \left( t,\overline{x}_{t}^{\varepsilon },%
\overline{u}_{t}^{\varepsilon }\right) \right) \right\} dt.%
\end{array}
\tag{3.12}  \label{3.12}
\end{equation}%
New we are to derive an estimate for the term similar to the right hand side
of the abov inquality with all the $\left( \overline{x}_{t}^{\varepsilon
},\left( \overline{u}_{t}^{\varepsilon },\overline{\eta }_{t}^{\varepsilon
}\right) \right) $ etc, replacing by $\left( x_{t}^{\varepsilon },\left(
u_{t}^{\varepsilon },\eta _{t}^{\varepsilon }\right) \right) $ etc. To this
end, we use similar method as in Zhou (\cite{zhoo}) we obtian the following
estimates:%
\begin{equation}
\begin{array}{c}
\mathbb{E}\int_{s}^{T}\left[ \overline{K}_{t}^{\varepsilon }\left( \sigma
\left( t,\overline{x}_{t}^{\varepsilon },u\right) -\sigma \left( t,\overline{%
x}_{t}^{\varepsilon },\overline{u}_{t}^{\varepsilon }\right) \right) \right.
\smallskip \\ 
-\left. K_{t}^{\varepsilon }\left( \sigma \left( t,x_{t}^{\varepsilon
},u\right) -\sigma \left( t,x_{t}^{\varepsilon },u_{t}^{\varepsilon }\right)
\right) \right] dt\leq C\varepsilon ^{\delta },%
\end{array}
\tag{3.13}  \label{3.13}
\end{equation}%
and

\begin{equation}
\begin{array}{l}
\mathbb{E}\int_{s}^{T}\left\{ \frac{1}{2}\left( \sigma \left( t,\overline{x}%
_{t}^{\varepsilon },u\right) -\sigma \left( t,\overline{x}_{t}^{\varepsilon
},\overline{u}_{t}^{\varepsilon }\right) \right) ^{\ast }\overline{Q}%
_{t}^{\varepsilon }\left( \sigma \left( t,\overline{x}^{\varepsilon
},u\right) -\sigma \left( t,\overline{x}_{t}^{\varepsilon },\overline{u}%
_{t}^{\varepsilon }\right) \right) \right. \smallskip \\ 
-\frac{1}{2}\left( \sigma \left( t,x_{t}^{\varepsilon },u\right) -\sigma
\left( t,x_{t}^{\varepsilon },u_{t}^{\varepsilon }\right) \right) ^{\ast
}Q_{t}^{\varepsilon }\left( \sigma \left( t,x^{\varepsilon },u\right)
-\sigma \left( t,x_{t}^{\varepsilon },u_{t}^{\varepsilon }\right) \right)
\smallskip \\ 
+\left[ \overline{\Psi }_{t}^{\varepsilon }\left( f\left( t,\overline{x}%
_{t}^{\varepsilon },u\right) -f\left( t,\overline{x}_{t}^{\varepsilon },%
\overline{u}_{t}^{\varepsilon }\right) \right) -\Psi _{t}^{\varepsilon
}\left( f\left( t,x_{t}^{\varepsilon },u\right) -f\left(
t,x_{t}^{\varepsilon },u_{t}^{\varepsilon }\right) \right) \right] \smallskip
\\ 
\left. +\left[ \ell \left( t,\overline{x}_{t}^{\varepsilon },u\right) -\ell
\left( t,\overline{x}_{t}^{\varepsilon },\overline{u}_{t}^{\varepsilon
}\right) \right] -\left[ \ell \left( t,x_{t}^{\varepsilon },u\right) -\ell
\left( t,x_{t}^{\varepsilon },u_{t}^{\varepsilon }\right) \right] \right\}
dt\smallskip \\ 
\leq C\varepsilon ^{\delta },%
\end{array}
\tag{3.14}  \label{3.14}
\end{equation}%
where $\left( \overline{\Psi }^{\varepsilon },\overline{K}^{\varepsilon
}\right) $\ and $\left( \overline{Q}^{\varepsilon },\overline{R}%
^{\varepsilon }\right) $\ are the solutions of adjoint equations (\ref{2.5})
and (\ref{2.6}) respectively corresponding to $\left( \overline{x}%
^{\varepsilon },\left( \overline{u}^{\varepsilon },\overline{\eta }%
^{\varepsilon }\right) \right) .$ The first variational inequality (\ref{3.6}%
) follows from combining, (\ref{3.12}) (\ref{3.13}) and.(\ref{3.14})

\noindent \textbf{Corollary 1.} Under the assumptions of Theorem 1, we have%
\begin{equation}
\mathbb{E}\int_{s}^{T}\mathcal{H}^{\left( x^{\varepsilon },u^{\varepsilon
}\right) }(t,x_{t}^{\varepsilon },u_{t}^{\varepsilon })dt\geq \sup_{u\in 
\mathbb{U}\left( \left[ s,T\right] \right) }\mathbb{E}\int_{s}^{T}\mathcal{H}%
^{\left( x^{\varepsilon },u^{\varepsilon }\right) }(t,x_{t}^{\varepsilon
},u_{t})dt-C\varepsilon ^{\delta }.  \tag{3.15}  \label{3.15}
\end{equation}

\noindent \textbf{Second variational inequality: }To obtain the second
variational inequality,\textit{\ }we define the following convex
perturbation $(\overline{u}_{t}^{\varepsilon },\overline{\eta }%
_{t}^{\varepsilon ,\theta })\in \mathbb{U}_{1}\times \mathbb{U}_{2}:$

\begin{equation}
(\overline{u}_{t}^{\varepsilon },\overline{\eta }_{t}^{\varepsilon ,\theta
})=\left( \overline{u}_{t}^{\varepsilon },\overline{\eta }_{t}^{\varepsilon
}+\theta \left( \xi -\overline{\eta }_{t}^{\varepsilon }\right) \right) , 
\tag{3.16}  \label{3.16}
\end{equation}%
where $\xi $ is an arbitrary element of the set $\mathbb{U}_{2}$. Using the
optimality of $\left( \overline{u}^{\varepsilon },\overline{\eta }%
^{\varepsilon }\right) $ to the new cost, $J^{\varepsilon }$ we have%
\begin{equation}
J^{\varepsilon }(s,y,\overline{u}^{\varepsilon },\overline{\eta }%
^{\varepsilon })\leq J^{\varepsilon }(s,y,\overline{u}_{t}^{\varepsilon },%
\overline{\eta }_{t}^{\varepsilon ,\theta }),  \tag{3.17}  \label{3.17}
\end{equation}%
a simple computation on $d_{2}(\overline{u}_{t}^{\varepsilon },\overline{%
\eta }_{t}^{\varepsilon ,\theta })$ we obtain%
\begin{equation*}
J(s,y,\overline{u}_{t}^{\varepsilon },\overline{\eta }_{t}^{\varepsilon
,\theta })-J\left( s,y,\overline{u}^{\varepsilon },\overline{\eta }%
^{\varepsilon }\right) \geq -C\theta \varepsilon ^{\frac{1}{3}}\geq -C\theta
\varepsilon ^{\delta }.
\end{equation*}%
Finally, arguing as in (\cite{seid}) for the left-hand side of the above
inequality, then we have%
\begin{equation*}
\lim_{\theta \rightarrow 0}\frac{1}{\theta }\left[ J(s,y,\overline{u}%
_{t}^{\varepsilon },\overline{\eta }_{t}^{\varepsilon ,\theta })-J\left( s,y,%
\overline{u}^{\varepsilon },\overline{\eta }^{\varepsilon }\right) \right] =%
\mathbb{E}\int_{s}^{T}(k_{t}+G_{t}^{\ast }\overline{\Psi }_{t}^{\varepsilon
})d\left( \eta -\overline{\eta }^{\varepsilon }\right) _{t},
\end{equation*}%
then the near-singular maximum condition follows%
\begin{equation}
-C\varepsilon ^{\delta }\leq \mathbb{E}\left[ \int_{s}^{T}(k_{t}+G_{t}^{\tau
}\overline{\Psi }_{t}^{\varepsilon })d\left( \eta -\overline{\eta }%
^{\varepsilon }\right) _{t}\right] .  \tag{3.18}  \label{3.18}
\end{equation}%
Now, we are to drive an estimate for the term similar to the right hand side
of (\ref{3.18}) with all $\left( \overline{x}^{\varepsilon },\left( 
\overline{u}^{\varepsilon },\overline{\eta }^{\varepsilon }\right) \right) $
ect., replacing by $\left( x^{\varepsilon },\left( u^{\varepsilon },\eta
^{\varepsilon }\right) \right) $ ect. We first estimate the following
difference:

\begin{eqnarray*}
&&\mathbb{E}\int_{s}^{T}(k_{t}+G_{t}^{\ast }\overline{\Psi }%
_{t}^{\varepsilon })d(\eta -\overline{\eta }^{\varepsilon })_{t}-\mathbb{E}%
\int_{s}^{T}(k_{t}+G_{t}^{\ast }\Psi _{t}^{\varepsilon })d(\eta -\eta
^{\varepsilon })_{t} \\
&=&\mathbb{E}\int_{s}^{T}G_{t}^{\ast }\left( \overline{\Psi }%
_{t}^{\varepsilon }-\Psi _{t}^{\varepsilon }\right) d\eta _{t}+\mathbb{E}%
\int_{s}^{T}(k_{t}+G_{t}^{\ast }\Psi _{t}^{\varepsilon })d\eta
_{t}^{\varepsilon } \\
&&-\mathbb{E}\int_{s}^{T}(k_{t}+G_{t}^{\ast }\overline{\Psi }%
_{t}^{\varepsilon })d\overline{\eta }_{t}^{\varepsilon }. \\
&=&\mathbb{E}\int_{s}^{T}G_{t}^{\ast }\left( \overline{\Psi }%
_{t}^{\varepsilon }-\Psi _{t}^{\varepsilon }\right) d\left( \eta -\overline{%
\eta }^{\varepsilon }\right) _{t}+\mathbb{E}\int_{s}^{T}(k_{t}+G_{t}^{\ast
}\Psi _{t}^{\varepsilon })d\left( \eta ^{\varepsilon }-\overline{\eta }%
^{\varepsilon }\right) _{t}.
\end{eqnarray*}%
Using the boundness of $G_{t},$ $k_{t}$, Lemma3, definition 1 ($\eta
_{s}=\eta _{s}^{\varepsilon }=\overline{\eta }_{s}^{\varepsilon }=0,$ $%
\mathbb{E}\left\vert \eta _{T}-\overline{\eta }_{T}^{\varepsilon
}\right\vert ^{2}+\mathbb{E}\left\vert \eta _{T}^{\varepsilon }-\overline{%
\eta }_{T}^{\varepsilon }\right\vert ^{2}<\infty )$ and the fact that $%
\mathbb{E}\left( \sup_{s\leq t\leq T}\left\vert \Psi _{t}^{\varepsilon
}\right\vert ^{2}\right) <C$ we have 
\begin{equation}
\mathbb{E}\int_{s}^{T}(k_{t}+G_{t}^{\ast }\overline{\Psi }_{t}^{\varepsilon
})d(\eta -\overline{\eta }^{\varepsilon })_{t}-\mathbb{E}%
\int_{s}^{T}(k_{t}+G_{t}^{\ast }\Psi _{t}^{\varepsilon })d(\eta -\eta
^{\varepsilon })_{t}\leq C\varepsilon ^{\delta }.  \tag{3.19}  \label{3.19}
\end{equation}%
Combining (\ref{3.18}) and (\ref{3.19}) the proof of inequality (\ref{3.7})
is complete.

\noindent \textbf{Corollary 2.} Under the assumptions of Theorem 1, we have

\begin{equation}
\mathbb{E}\int_{s}^{T}(k_{t}+G_{t}^{\ast }\Psi _{t}^{\varepsilon })d\eta
_{t}^{\varepsilon }\leq \inf_{\eta \in \mathbb{U}_{2}\left( \left[ s,T\right]
\right) }\mathbb{E}\int_{s}^{T}(k_{t}+G_{t}^{\ast }\Psi _{t}^{\varepsilon
})d\eta _{t}+C\varepsilon ^{\delta }.  \tag{3.20}  \label{3.20}
\end{equation}

\section{Sufficient near-optimality conditions}

In this section, we will prove that under an additional assumptions, the
near-maximum condition on the Hamiltonian function is sufficient for
near-optimality. We assume:

\begin{description}
\item[(H4)] $\rho $ is differentiable in $u$ for $\varphi =f,\sigma ,\ell $
and there is a constante $C$ such that%
\begin{equation}
\left\vert \varphi (t,x,u,)-\varphi (t,x,u^{\prime },)\right\vert
+\left\vert \varphi _{u}(t,x,u,)-\varphi _{u}(t,x,u^{\prime },)\right\vert
\leq C\left\vert u-u^{\prime }\right\vert .  \tag{4.1}  \label{4.1}
\end{equation}
\end{description}

\noindent \textbf{Theorem 2. }\textit{Assume the }$H\left( t,\cdot ,\cdot
,\Psi _{t}^{\varepsilon },K_{t}^{\varepsilon }\right) $\textit{\ is concave
for }$a.e.$\textit{\ }$t\in \left[ s,T\right] ,$\textit{\ }$\mathbb{P}-a.s$%
\textit{, and }$h$\textit{\ is convex. Let }$\left( \Psi _{t}^{\varepsilon
},K_{t}^{\varepsilon }\right) ,$\textit{\ }$\left( Q_{t}^{\varepsilon
},R_{t}^{\varepsilon }\right) $\textit{\ be the solution of the adjoint
equation (\ref{2.5})-(\ref{2.6}) associated with }$\left( u^{\varepsilon
},\eta ^{\varepsilon }\right) .$\textit{\ If for some }$\varepsilon >0$%
\textit{\ and for any }$\left( u,\eta \right) \in \mathbb{U}:$%
\begin{equation}
\mathbb{E}\int_{s}^{T}\mathcal{H}^{\left( x^{\varepsilon },u^{\varepsilon
}\right) }(t,x_{t}^{\varepsilon },u_{t}^{\varepsilon })dt\geq \sup_{u\in 
\mathbb{U}_{1}\left( \left[ s,T\right] \right) }\mathbb{E}\int_{s}^{T}%
\mathcal{H}^{\left( x^{\varepsilon },u^{\varepsilon }\right)
}(t,x_{t}^{\varepsilon },u_{t})dt-\varepsilon ,  \tag{4.2}  \label{4.2}
\end{equation}%
and%
\begin{equation}
\mathbb{E}\left[ \int_{s}^{T}k_{t}d\left( \eta -\eta ^{\varepsilon }\right)
_{t}\right] \geq -C\varepsilon ^{\frac{1}{2}},  \tag{4.3}  \label{4.3}
\end{equation}%
\textit{then we have}%
\begin{equation}
J\left( s,y,u^{\varepsilon },\eta ^{\varepsilon }\right) \leq \inf_{\left(
u,\eta \right) \in \mathbb{U}\left( \left[ s,T\right] \right) }J\left(
s,y,u,\eta \right) +C\varepsilon ^{\frac{1}{2}},  \tag{4.4}  \label{4.4}
\end{equation}%
\textit{where }$C$\textit{\ is a positive constant independent of }$%
\varepsilon .$

\noindent \textbf{Proof. }First, define the cost functional 
\begin{equation}
J\left( s,y,u,\eta \right) =J_{1}\left( s,y,u\right) +J_{2}\left( s,\eta
\right) ,  \tag{4.5}  \label{4.5}
\end{equation}%
where%
\begin{equation*}
J_{1}\left( s,y,u\right) =\mathbb{E}\left[ h\left( x_{T}\right)
+\int_{s}^{T}\ell \left( t,x_{t},u_{t}\right) dt\right] ,
\end{equation*}%
and%
\begin{equation*}
J_{2}\left( s,\eta \right) =\mathbb{E}\left[ \int_{s}^{T}k_{t}d\eta _{t}%
\right] .
\end{equation*}%
Let us fix $\varepsilon >0,$ Define a new metric $\widetilde{d}$ on $\mathbb{%
U}\left( \left[ s,T\right] \right) $\ as follows: for any $\left( u,\eta
\right) \ $and $\left( v,\xi \right) \in \mathbb{U}:$%
\begin{equation*}
\widetilde{d}\left( \left( u,\eta \right) ,\left( v,\xi \right) \right) =%
\widetilde{d_{1}}\left( u,v\right) +d_{2}\left( \eta ,\xi \right) ,
\end{equation*}%
where%
\begin{equation}
\widetilde{d_{1}}\left( u,v\right) =\mathbb{E}\left[ \int_{s}^{T}\varsigma
_{t}^{\varepsilon }\left\vert u_{t}-v_{t}\right\vert dt\right] ,  \tag{4.6}
\label{4.6}
\end{equation}%
and%
\begin{equation}
\varsigma _{t}^{\varepsilon }=1+\left\vert \Psi _{t}^{\varepsilon
}\right\vert +\left\vert K_{t}^{\varepsilon }\right\vert +\left\vert
Q_{t}^{\varepsilon }\right\vert +\left\vert R_{t}^{\varepsilon }\right\vert
\geq 1.  \tag{4.7}  \label{4.7}
\end{equation}%
Obviously $\widetilde{d_{1}}$ is a metric on $(\mathbb{U}_{1},\widetilde{d}%
_{1})$, and it is a complete metric as a weighted $\mathbb{L}^{1}$ norm.
Hence $(\mathbb{U},\widetilde{d})$ as a product of two complete metric
spaces is a complete metric space under $\widetilde{d}$.

\noindent Define a functional $\Upsilon $ on $\mathbb{U}_{1}\left( \left[ s,T%
\right] \right) $ by%
\begin{eqnarray*}
\Upsilon \left( u\right) &=&\mathbb{E}\int_{s}^{T}H\left( t,x^{\varepsilon
},u,\Psi ^{\varepsilon },K^{\varepsilon },Q^{\varepsilon }\right) dt \\
&=&\mathbb{E}\int_{s}^{T}\mathcal{H}^{\left( x^{\varepsilon },u^{\varepsilon
}\right) }\left( t,x^{\varepsilon },u\right) dt,
\end{eqnarray*}%
a simple computation shows that%
\begin{equation*}
\left\vert \Upsilon \left( u\right) -\Upsilon \left( v\right) \right\vert
\leq C\mathbb{E}\left[ \int_{s}^{T}\varsigma _{t}^{\varepsilon }\left\vert
u_{t}-v_{t}\right\vert dt\right] ,
\end{equation*}%
which implies that $\Upsilon $ is continuous on $\mathbb{U}_{1}\left( \left[
s,T\right] \right) $ with respct to $\widetilde{d_{1}}.$ Now by using (\ref%
{4.2}) and Ekeland variational principle, there exists a $\overline{u}%
^{\varepsilon }\in \mathbb{U}_{1}\left( \left[ s,T\right] \right) $ such that

\begin{equation*}
\widetilde{d_{1}}(\overline{u}^{\varepsilon ,},u^{\varepsilon })\leq
\varepsilon ^{\frac{1}{2}},
\end{equation*}%
and%
\begin{equation}
\mathbb{E}\int_{s}^{T}\widetilde{\mathcal{H}}(t,x_{t}^{\varepsilon },%
\overline{u}_{t}^{\varepsilon })dt=\max_{u\in \mathbb{U}_{1}\left( \left[ s,T%
\right] \right) }\mathbb{E}\int_{s}^{T}\widetilde{\mathcal{H}}%
(t,x_{t}^{\varepsilon },u_{t})dt,  \tag{4.8}  \label{4.8}
\end{equation}%
where%
\begin{equation}
\widetilde{\mathcal{H}}(t,x,u)=\mathcal{H}^{\left( x^{\varepsilon
},u^{\varepsilon }\right) }(t,x,u)-\varepsilon ^{\frac{1}{2}}\varsigma
_{t}^{\varepsilon }\left\vert u-\overline{u}_{t}^{\varepsilon }\right\vert .
\tag{4.9}  \label{4.9}
\end{equation}%
The maximum condition (\ref{4.8}) implies a\ pointwise maximum condition
namely, for $a.e.$ $t\in \left[ s,T\right] $ and $\mathbb{P}-a.s,$%
\begin{equation*}
\widetilde{\mathcal{H}}(t,x_{t}^{\varepsilon },\overline{u}_{t}^{\varepsilon
})=\max_{u\in \mathbb{A}}\widetilde{\mathcal{H}}(t,x^{\varepsilon },u).
\end{equation*}%
Using Proposition A1 (Appendix), then we have%
\begin{equation*}
\partial _{u}\widetilde{\mathcal{H}}(t,x_{t}^{\varepsilon },\overline{u}%
_{t}^{\varepsilon })\ni 0.
\end{equation*}%
Since $\left\vert u-\overline{u}_{t}^{\varepsilon }\right\vert $ is not
differentiable in $\overline{u}_{t}^{\varepsilon }$ (locally Lipschitz),
then we use Proposition A1 (Appendix) we get

\begin{equation*}
\partial _{u}\left\{ \varepsilon ^{\frac{1}{2}}\varsigma _{t}^{\varepsilon
}\left\vert u-\overline{u}_{t}^{\varepsilon }\right\vert \right\} =\left[
-\varepsilon ^{\frac{1}{2}}\varsigma _{t}^{\varepsilon },\text{ }\varepsilon
^{\frac{1}{2}}\varsigma _{t}^{\varepsilon }\right] .
\end{equation*}%
By using (\ref{4.9}) and fact that the Clarke's generalized gradient of the
sum of two functions is contained in the sum of the Clarke's generalized
gradient of the two functions, we get

\begin{equation*}
\partial _{u}\widetilde{\mathcal{H}}(t,x_{t}^{\varepsilon },\overline{u}%
_{t}^{\varepsilon })\subset \partial _{u}\mathcal{H}^{\left( x^{\varepsilon
},u^{\varepsilon }\right) }(t,x_{t}^{\varepsilon },\overline{u}%
_{t}^{\varepsilon })+\left[ -\varepsilon ^{\frac{1}{2}}\varsigma
_{t}^{\varepsilon },\text{ }\varepsilon ^{\frac{1}{2}}\varsigma
_{t}^{\varepsilon }\right] .
\end{equation*}%
Applying the similar method as in (\cite{zhoo}) for the rest of the proof,
we obtian, for an arbitrary $u$

\begin{equation}
J_{1}\left( s,y,u^{\varepsilon }\right) \leq J_{1}\left( s,y,u\right)
+C\varepsilon ^{\frac{1}{2}}.  \tag{4.10}  \label{4.10}
\end{equation}%
Now, by using (\ref{4.3}) we get%
\begin{eqnarray*}
J_{2}\left( s,\eta ^{\varepsilon }\right) &=&\mathbb{E}\left[
\int_{s}^{T}k_{t}d\eta _{t}^{\varepsilon }\right] \\
&\leq &\mathbb{E}\left[ \int_{s}^{T}k_{t}d\eta _{t}\right] +C\varepsilon ^{%
\frac{1}{2}},
\end{eqnarray*}%
which implies that for an arbitrary $\eta $%
\begin{equation}
J_{2}\left( s,\eta ^{\varepsilon }\right) \leq J_{2}\left( s,\eta \right)
+C\varepsilon ^{\frac{1}{2}}.  \tag{4.11}  \label{4.11}
\end{equation}%
Combining (\ref{4.10}), (\ref{4.11}) and (\ref{4.5}) we arrive at%
\begin{equation*}
J\left( s,y,u^{\varepsilon },\eta ^{\varepsilon }\right) \leq J\left(
s,y,u,\eta \right) +C\varepsilon ^{\frac{1}{2}}.
\end{equation*}%
Since $\left( u,\eta \right) $ is arbitrary, the desired result follows.

\noindent \textbf{Corollary 3.\ }Let the assumptions of Theorem 2 holds. A
sufficient conditions for an admissible pair $(x^{\varepsilon
},u^{\varepsilon },\eta ^{\varepsilon })$ to be $\varepsilon -$optimal is%
\begin{eqnarray*}
&&\mathbb{E}\left\{ \int_{s}^{T}\mathcal{H}^{\left( x^{\varepsilon
},u^{\varepsilon }\right) }(t,x_{t}^{\varepsilon },u_{t}^{\varepsilon
})dt-\int_{s}^{T}k_{t}d\eta _{t}^{\varepsilon }\right\} \\
&\geq &\sup_{\left( u,\eta \right) \in \mathbb{U}\left( \left[ s,T\right]
\right) }\mathbb{E}\left\{ \int_{s}^{T}\mathcal{H}^{\left( x^{\varepsilon
},u^{\varepsilon }\right) }(t,x_{t}^{\varepsilon
},u_{t})dt-\int_{s}^{T}k_{t}d\eta _{t}\right\} -\left( \frac{\varepsilon }{C}%
\right) ^{2}.
\end{eqnarray*}%
\noindent \textbf{Example 1.} Consider the one-dimensional stochastic
control problem: $n=l=1,$ $G_{t}=1,$ $\mathbb{A}_{1}=\left[ 0,1\right] $, $%
\mathbb{A}_{2}=\left[ 0,1\right] ,$ $\eta _{1}=1,$%
\begin{equation}
\left\{ 
\begin{array}{c}
dx_{t}=u_{t}dW_{t}+d\eta _{t}, \\ 
\multicolumn{1}{l}{x_{0}=0,}%
\end{array}%
\right.  \tag{4.12}  \label{4.12}
\end{equation}%
and the cost functional being%
\begin{equation}
J\left( s,y,u,\eta \right) =\mathbb{E}\left[ \frac{1}{2}x_{1}^{2}-%
\int_{0}^{1}u_{t}dt+\int_{0}^{1}k_{t}d\eta _{t}\right] ,  \tag{4.13}
\label{4.13}
\end{equation}%
For a given admissible pair $\left( x^{\varepsilon },(u^{\varepsilon },\eta
^{\varepsilon })\right) $, the correspending second-order adjoint equation is%
\begin{equation}
\left\{ 
\begin{array}{c}
dQ_{t}^{\varepsilon }=R_{t}^{\varepsilon }dW_{t} \\ 
\multicolumn{1}{l}{Q_{1}^{\varepsilon }=1,}%
\end{array}%
\right.  \tag{4.14}  \label{4.14}
\end{equation}%
By the uniqueness of this solution, $\left( Q^{\varepsilon },R^{\varepsilon
}\right) =\left( 1,0\right) $, then for any admissible control $\left(
u,\eta \right) $ we have%
\begin{eqnarray*}
\mathcal{H}^{\left( x^{\varepsilon },u^{\varepsilon }\right)
}(t,x_{t}^{\varepsilon },u_{t}) &=&u_{t}-\left( K_{t}^{\varepsilon
}-Q_{t}^{\varepsilon }u_{t}^{\varepsilon }\right) u_{t}-\tfrac{1}{2}%
Q_{t}^{\varepsilon }u_{t}^{2} \\
&=&\tfrac{1}{2}\left[ \left( u_{t}^{\varepsilon }-K_{t}^{\varepsilon
}+1\right) ^{2}-\left( u_{t}-u_{t}^{\varepsilon }+K_{t}^{\varepsilon
}-1\right) ^{2}\right] .
\end{eqnarray*}%
Replacing $u_{t}=u_{t}^{\varepsilon }$, we get%
\begin{eqnarray*}
\mathcal{H}^{\left( x^{\varepsilon },u^{\varepsilon }\right)
}(t,x_{t}^{\varepsilon },u_{t}^{\varepsilon }) &=&\tfrac{1}{2}\left[ \left(
u_{t}^{\varepsilon }-\left( K_{t}^{\varepsilon }-1\right) \right)
^{2}-\left( K_{t}^{\varepsilon }-1\right) ^{2}\right] \\
&=&\tfrac{1}{2}\left[ \left( u_{t}^{\varepsilon }\right)
^{2}+2u_{t}^{\varepsilon }\left( 1-K_{t}^{\varepsilon }\right) \right] .
\end{eqnarray*}%
Hence by simple computation shows that if 
\begin{equation}
u_{t}^{\varepsilon }-\left( K_{t}^{\varepsilon }-1\right) \in \left[ 0,1%
\right] ,  \tag{4.15}  \label{4.15}
\end{equation}%
then\ (\ref{3.15}) and (\ref{3.20}) gives%
\begin{eqnarray*}
&&\mathbb{E}\int_{0}^{1}\tfrac{1}{2}\left[ \left( u_{t}^{\varepsilon
}\right) ^{2}+2u_{t}^{\varepsilon }\left( 1-K_{t}^{\varepsilon }\right) %
\right] dt \\
&\geq &\sup_{u\in \left[ 0,1\right] }\mathbb{E}\int_{0}^{1}\tfrac{1}{2}\left[
\left( u_{t}^{\varepsilon }-K_{t}^{\varepsilon }+1\right) ^{2}-\left(
u_{t}-u_{t}^{\varepsilon }+K_{t}^{\varepsilon }-1\right) ^{2}\right]
dt-C\varepsilon ^{\delta },
\end{eqnarray*}%
and%
\begin{equation*}
\mathbb{E}\int_{s}^{T}(k_{t}+\Psi _{t}^{\varepsilon })d\eta
_{t}^{\varepsilon }\leq \inf_{\eta \in \mathbb{U}_{2}\left( \left[ s,T\right]
\right) }\mathbb{E}\int_{s}^{T}(k_{t}+\Psi _{t}^{\varepsilon })d\eta
_{t}+C\varepsilon ^{\delta },
\end{equation*}%
thus%
\begin{equation}
\mathbb{E}\int_{0}^{1}\left( K_{t}^{\varepsilon }-1\right) ^{2}dt\leq
C\varepsilon ^{\delta }.  \tag{4.16}  \label{4.16}
\end{equation}%
We denote $\mathfrak{B}=\left\{ \left( w,t\right) \in \Omega \times \left[
0,1\right] :k_{t}+\Psi _{t}^{\varepsilon }\geq -C\varepsilon ^{\frac{1}{2}%
}\right\} $ and define%
\begin{equation*}
d\eta _{t}=\left\{ 
\begin{array}{l}
0\text{ if }k_{t}+\Psi _{t}^{\varepsilon }\geq -C\varepsilon ^{\frac{1}{2}%
}\smallskip \\ 
d\eta _{t}^{\varepsilon }\text{ otherwise.}%
\end{array}%
\right.
\end{equation*}%
Hence a simple calculation shows that%
\begin{equation*}
\mathbb{E}\int_{0}^{1}(k_{t}+\Psi _{t}^{\varepsilon })d\left( \eta -\eta
^{\varepsilon }\right) _{t}=\mathbb{E}\int_{0}^{1}(k_{t}+\Psi
_{t}^{\varepsilon })\mathbf{\chi }_{\mathfrak{B}}d\left( -\eta ^{\varepsilon
}\right) _{t}\geq -C\varepsilon ^{\frac{1}{2}},
\end{equation*}%
which implies that%
\begin{equation}
\mathbb{E}\int_{0}^{1}(k_{t}+\Psi _{t}^{\varepsilon })\mathbf{\chi }_{%
\mathfrak{B}}d\eta _{t}^{\varepsilon }\leq C\varepsilon ^{\delta } 
\tag{4.17}  \label{4.17}
\end{equation}%
It worth montioning that the above conditions reveals the minimum
qualification for the pair $\left( x^{\varepsilon },(u^{\varepsilon },\eta
^{\varepsilon })\right) $ to be $\varepsilon -$optimal. As an example, the
admissible controls $(u^{\varepsilon },\eta ^{\varepsilon })=(1-\varepsilon
^{\frac{1}{2}},\eta ^{\varepsilon })$ are condidate $\varepsilon -$%
optimality, where $\varepsilon >0$ is sufficiently small. Note that the
first-order adjoint equation is%
\begin{equation}
\left\{ 
\begin{array}{c}
d\Psi _{t}^{\varepsilon }=K_{t}^{\varepsilon }dW_{t} \\ 
\multicolumn{1}{l}{\Psi _{1}^{\varepsilon }=x_{1}^{\varepsilon },}%
\end{array}%
\right.  \tag{4.18}  \label{4.18}
\end{equation}%
with the corresponding trajectories $x_{t}^{\varepsilon }=(1-\varepsilon ^{%
\frac{1}{2}})W_{t}+\eta _{t}^{\varepsilon }$ then the unique solution pair
of the first-order adjoint equation will be $\left( \Psi _{t}^{\varepsilon
},K_{t}^{\varepsilon }\right) =\left( (1-\varepsilon ^{\frac{1}{2}%
})W_{t}+1,(1-\varepsilon ^{\frac{1}{2}})\right) .$ Hence (\ref{4.15}) and (%
\ref{4.16}) will be satisfied.

\noindent Conversely, for the sufficient part, since the hamiltonian $%
H\left( t,x,u,p,q\right) =u-qu$ is concave in $\left( x,u\right) $, we use
Theorem 1 to conclude that $u_{t}^{\varepsilon }=(1-\left( \frac{\varepsilon 
}{C}\right) ^{2})$ is a condidate to be $\varepsilon -$optimality for
sufficiently small $\varepsilon $ is indeed an $\varepsilon -$optimal
control.

\begin{description}
\item[{\protect\LARGE Appendix}] 
\end{description}

The following result gives some basic properties of the Clarke's gneralized
gradient.

\noindent \textbf{Proposition A1. }If $f:\mathbb{R}^{n}\rightarrow \mathbb{R}
$ is locally Lipschitz at $x\in \mathbb{R}^{n}$, then the following
statements holds

\noindent (1) $\partial _{x}f\left( x\right) $ is nonempty, compact, and
convex set in $\mathbb{R}^{n}$

\noindent (2) $\partial _{x}\left( -f\right) \left( x\right) =-\partial
_{x}\left( f\right) \left( x\right) ,$

\noindent (3) $0\in \partial _{x}\left( f\right) \left( x\right) $ if $f$
attains a local minimum or maximum at $x$

\noindent (4) If $f$\ is Fr\^{e}chet-differentiable at $x$, then $\partial
_{x}f\left( x\right) =\left\{ f^{\prime }\left( x\right) \right\} .$

\noindent (5) If $f,$ $g:\mathbb{R}^{n}\rightarrow \mathbb{R}$ are locally
Lipschitz at $x\in \mathbb{R}^{d}$, then%
\begin{equation*}
\partial _{x}\left( f+g\right) \left( x\right) \subset \partial _{x}f\left(
x\right) +\partial _{x}g\left( x\right) .
\end{equation*}

\noindent See (\cite{yong} Lemma 2.3), \cite{clarke}\ for the detailed proof
of the above Proposition.

\noindent As an example, the Clarke's generalized gradient of the absolute
value function $f:x\mapsto \left\vert x\right\vert $ which is continuously
differentiable everywhere except at $0$. Since $f^{\prime }\left( x\right)
=1 $ for $x>0$ and, $f^{\prime }\left( x\right) =-1$ for $x<0$ then the
Clarke's generalized gradient of $f$ at $x=0$ is given by%
\begin{equation*}
\partial _{x}f\left( 0\right) =\overline{co}\left\{ -1,1\right\} =\left[ -1,1%
\right] ,
\end{equation*}

\end{document}